\documentclass[11pt, oneside]{amsart}

\usepackage[utf8]{inputenc}

\usepackage{amsfonts, amstext, amsmath, amsthm, amscd, amssymb}
\usepackage{mathrsfs}
\usepackage{tikz}
\usetikzlibrary{knots}
\tikzset{point/.style={insert path={ node[scale=2.5*sqrt(\pgflinewidth)]{.} }}}
\usepackage{url}
\usepackage[all]{xypic}

\usepackage{color}
\usepackage[colorlinks,
    linkcolor={red!50!black},
    citecolor={blue!50!black},
    urlcolor={blue!80!black}]{hyperref}

\usepackage[margin=1.3in]{geometry}

\usepackage{enumerate}

\renewcommand{\setminus}{{\smallsetminus}}

\newcommand{\bp}{\begin{pmatrix}}
\newcommand{\ep}{\end{pmatrix}}
\newcommand{\be}{\begin{equation}}
\newcommand{\ee}{\end{equation}}

\numberwithin{equation}{section}

\theoremstyle{plain}
\newtheorem{theorem}[equation]{Theorem}
\newtheorem{lemma}[equation]{Lemma}

\newtheorem{conjecture}[equation]{Conjecture}

\theoremstyle{definition}

\newtheorem{definition}[equation]{Definition}
\newtheorem{problem}[equation]{Problem}

\newtheorem{remark}[equation]{Remark}

\theoremstyle{remark}

\numberwithin{equation}{section}

\def\Z{\mathbb Z}

\def\sm{\setminus}

\def\bp{\begin{pmatrix}}
\def\ep{\end{pmatrix}}
\def\ba{\begin{array}}
\def\ea{\end{array}}
\def\bn{\begin{enumerate}}
\def\en{\end{enumerate}}

\DeclareMathOperator\lk{lk}

\begin{document}

\title{The round handle problem}

\author{Min Hoon Kim}
\address{Center for Research in Topology, Department of Mathematics, POSTECH, Pohang 37673, Republic of Korea}
\email{kminhoon@gmail.com}

\author{Mark Powell}
\address{Department of Mathematical Sciences, Durham University, United Kingdom}
\email{mark.a.powell@durham.ac.uk}

\author{Peter Teichner}
\address{Max Planck Institut f\"{u}r Mathematik, Vivatsgasse 7, 53111 Bonn, Germany}
\email{teichner@mpim-bonn.mpg.de}


\def\subjclassname{\textup{2010} Mathematics Subject Classification}
\expandafter\let\csname subjclassname@1991\endcsname=\subjclassname
\expandafter\let\csname subjclassname@2000\endcsname=\subjclassname

\begin{abstract}
We present the Round Handle Problem (RHP), proposed by Freedman and Krushkal. It asks whether a collection of links, which contains the Generalised Borromean Rings (GBRs), are slice in a $4$-manifold $R$ constructed from adding round handles to the four ball.  A negative answer would contradict the union of the surgery conjecture and the $s$-cobordism conjecture for $4$-manifolds with free fundamental group.
\end{abstract}

\maketitle

\section{Statement of the RHP}

We give an alternative proof of the connection of the Round Handle Problem to the topological surgery and $s$-cobordism conjectures (these will all be recalled below).  The Round Handle Problem (RHP) was formulated in \cite[Section~5.1]{Freedman-Krushkal-2016-A}.    We give a shorter and easier argument that the above mentioned conjectures imply a positive answer to the RHP.

Let $L = L_1 \sqcup \cdots \sqcup L_m$ be an oriented ordered link in $S^3$ with vanishing pairwise linking numbers.  We will be particularly concerned with the Generalised Borromean Rings (GBRs). By definition these are the collection of links arising from iterated Bing doubling starting with a Hopf link. An example is shown in Figure~\ref{figure:GBR}.

\begin{figure}[htb]
\centering
\begin{tikzpicture} [scale=0.8,line width=1pt,double distance=1pt,join=round,
  over/.style={draw=white,double=black,double distance=1pt,line width=1.8pt},
  overred/.style={draw=white,double=red,double distance=1pt,line width=2pt},
  overblue/.style={draw=white,double=blue,double distance=1pt,line width=2pt},
  overwhite/.style={draw=white,double=white,double distance=1pt,line width=2pt},
  overparallels/.style 2 args={
    decoration={
      markings,
      mark=between positions 0 and 1-0.5*#2 step #2 with {
        \path +(0,#1) coordinate (ppp\pseqno) +(0,-#1) coordinate (qqq\pseqno);
        \xdef\lastno{\pseqno}
      },
      mark=at position .999 with {
        \pgfmathtruncatemacro{\pcnt}{\lastno+1}
        \fill
        +(0,#1) coordinate (ppp\pcnt) circle(.5pt) +(0,-#1) coordinate (qqq\pcnt) circle(.5pt);
        \gdef\ppps{}\gdef\qqqs{}
        \foreach \i in {1,...,\pcnt} {\xdef\ppps{\ppps (ppp\i)}\xdef\qqqs{\qqqs (qqq\i)}}
        \draw plot[smooth] coordinates {\ppps};
        \draw plot[smooth] coordinates {\qqqs};
        \fill (ppp1) circle(.5pt) (qqq1) circle(.5pt);
      },
    },
    postaction=decorate,
  },
  overparallels/.default={.6}{.05},
  ]

  \def\drawoverpiece{
    \path[overparallels] (A0) ..controls+(0,-3)and+(0,3).. (B0);
  }

  \def\whiteheadminus#1{
  \draw (0,2)\c[left]{AA#1} arc (90:180:2)--+(0,-2) arc (180:360:8 and 6)--+(0,6)arc (0:180:8 and 6)--+(0,-2)\c[left]{AB#1};
  \draw[over] (AB#1) arc (180:360:2)--+(0,2) arc(180:0:4 and 3)--+(0,-6) arc (360:180:4 and 3)--+(0,2) \c[left]{AC#1};
  \draw[over](AC#1) arc(0:90:2);
  }

  \def\whiteheadplus#1{
  \draw(2,0) arc(0:90:2);
  \draw (0,2)\c[left]{AA#1} arc (90:180:2)--+(0,-2) arc (180:360:8 and 6)--+(0,6)arc (0:180:8 and 6)--+(0,-2)\c[left]{AB#1};
  \draw[over] (AB#1) arc (180:360:2)--+(0,2) arc(180:0:4 and 3)--+(0,-6) arc (360:180:4 and 3)--+(0,2) \c[left]{AC#1};
  \draw[over] (0,2) arc(90:180:2);
  }

  \def\smallwhiteheadplus#1#2{
  \draw[over](2*#2,0) arc(0:90:2*#2);
  \draw[over] (0,2*#2)\c[left]{AA#1} arc (90:180:2*#2)--+(0,-2*#2) arc (180:360:8 *#2 and 6*#2)--+(0,6*#2)arc (0:180:8*#2 and 6*#2)--+(0,-2*#2)\c[left]{AB#1};
  \draw[over] (AB#1) arc (180:360:2*#2)--+(0,2*#2) arc(180:0:4*#2 and 3*#2)--+(0,-6*#2) arc (360:180:4*#2 and 3*#2)--+(0,2*#2) \c[left]{AC#1};
  \draw[over] (0,2*#2) arc(90:180:2*#2)--+(0,-#2);
  \filldraw[blue](14*#2,#2)  circle (1.3pt);
  }

  \def\smallwhiteheadminus#1#2{
  \draw[over](2*#2,0) arc(0:90:2*#2);
  \draw[over] (0,2*#2)\c[left]{AA#1} arc (90:180:2*#2)--+(0,-2*#2) arc (180:360:8 *#2 and 6*#2)--+(0,6*#2)arc (0:180:8*#2 and 6*#2)--+(0,-2*#2)\c[left]{AB#1};
  \draw[over] (AB#1) arc (180:360:2*#2)--+(0,2*#2) arc(180:0:4*#2 and 3*#2)--+(0,-6*#2) arc (360:180:4*#2 and 3*#2)--+(0,2*#2) \c[left]{AC#1};
  \draw[over] (0,2*#2) arc(90:0:2*#2)--+(0,-#2);
  \filldraw[blue](14*#2,#2)  circle (1.3pt);
  }
  \def\vinsidetorus#1{
\draw(0.07*#1  ,-0.4*#1  ) to [out=left,in=down] (-0.17*#1  ,0 )to [out=up,in=left] (0.07*#1  ,.4*#1  );
\draw(-0.09*#1  ,-0.47*#1   ) to [out=right,in=down] (0.17*#1  , 0 )to [out=up,in=right] (-0.09*#1  , .47*#1  );
  }
  \def\vinsidetorusred#1{
\draw[red](0.07*#1  ,-0.4*#1  ) to [out=left,in=down] (-0.17*#1  ,0 )to [out=up,in=left] (0.07*#1  ,.4*#1  );
\draw[red](-0.09*#1  ,-0.47*#1   ) to [out=right,in=down] (0.17*#1  , 0 )to [out=up,in=right] (-0.09*#1  , .47*#1  );
  }
  \def\hinsidetorus#1{
  \begin{scope}[rotate=-90]
\draw(0.07*#1  ,-0.4*#1  ) to [out=left,in=down] (-0.17*#1  ,0 )to [out=up,in=left] (0.07*#1  ,.4*#1  );
\draw(-0.09*#1  ,-0.47*#1   ) to [out=right,in=down] (0.17*#1  , 0 )to [out=up,in=right] (-0.09*#1  , .47*#1  );
\end{scope}
  }
  \tikzset{
      partial ellipse/.style args={#1:#2:#3}{
          insert path={+ (#1:#3) arc (#1:#2:#3)}
      }
  }

\begin{scope}[yshift=-10.5cm]

\draw(0.8,-1.4) arc (270:90:1.6) arc(270:450:0.1) arc(90:270:1.8) arc(270:450:0.1);
\draw[over] (0.9, -1.4) arc(90:270:0.1)--++(0.4,0) arc(270:450:0.4)--++(-0.5,0) arc(270:90:1)--++(0.5,0) arc(270:450:0.4)--++(-0.4,0) arc(90:270:0.1)--++(0.4,0) arc(450:270:0.2)--++(-0.5,0) arc(90:270:1.2)--++(0.5,0)arc (450:270:0.2)--++(-0.4,0);

\draw (1.1,1.6) arc (180:90:0.4)--+(0.4,0) arc(450:270:0.1)--+(-0.4,0) arc(90:270:0.2)--+(0.5,0) arc(450:270:1.2)--+(-0.5,0) arc(90:270:0.2)--+(0.5,0) arc(270:450:1.6);

\draw[over] (1.5,1.2) arc (270:90:0.4)--+(0.4,0) arc(450:270:0.1)--+(-0.4,0) arc(90:270:0.2)--+(0.5,0) arc(450:270:1.2)--+(-0.5,0) arc(90:270:0.2)--+(0.4,0) arc(450:270:0.1)--++(-0.4,0) arc(270:90:0.4)--++(0.5,0) arc(270:450:1)--++(-0.5,0);
\draw[over](2,-1.4) arc (270:450:1.6) arc(270:90:0.1) arc(450:270:1.8) arc(270:90:0.1);

\draw[over] (1.5,1.6)arc (360:270:0.2) --+(-0.4,0);
\draw[over] (1.7,1.6)arc (360:270:0.4) --+(-0.4,0);
\draw[over] (1.5,-1.2)arc (0:90:0.2) --+(-0.4,0);
\draw[over] (1.7,-1.2)arc (0:90:0.4) --+(-0.4,0);
\draw[over] (2,1.9)arc (0:90:0.1) --+(-0.2,0);
\draw[over] (0.9,1.9)arc (0:90:0.1) arc(90:100:1.8);
\draw[over] (0.9,-1.5)arc (360:270:0.1) arc(270:260:1.8);
\draw[over] (1.9,1.9)arc (180:270:0.1) arc(450:360:1.6);

\draw[over] (2,-1.5)arc (360:270:0.1) --+(-0.2,0);
\end{scope}

\begin{scope}[yshift=-10.5cm,xshift=-2.8cm]

\draw(0.8,-1.4) arc (270:90:1.6) arc(270:450:0.1) arc(90:270:1.8) arc(270:450:0.1);

\draw[over] (0.9, -1.4) arc(90:270:0.1)--++(0.4,0) arc(270:450:0.4)--++(-0.5,0) arc(270:90:1)--++(0.5,0) arc(270:450:0.4)--++(-0.4,0) arc(90:270:0.1)--++(0.4,0) arc(450:270:0.2)--++(-0.5,0) arc(90:270:1.2)--++(0.5,0)arc (450:270:0.2)--++(-0.4,0);
\draw[blue](2,-1.4) arc (270:450:1.6) arc(270:90:0.1) arc(450:270:1.8) arc(270:90:0.1);

\draw (1.1,1.6) arc (180:90:0.4)--+(0.4,0) arc(450:270:0.1)--+(-0.4,0) arc(90:270:0.2)--+(0.5,0) arc(450:270:1.2)--+(-0.5,0) arc(90:270:0.2)--+(0.5,0) arc(270:450:1.6);
\draw[over] (1.5,1.2) arc (270:90:0.4)--+(0.4,0) arc(450:270:0.1)--+(-0.4,0) arc(90:270:0.2)--+(0.5,0) arc(450:270:1.2)--+(-0.5,0) arc(90:270:0.2)--+(0.4,0) arc(450:270:0.1)--++(-0.4,0) arc(270:90:0.4)--++(0.5,0) arc(270:450:1)--++(-0.5,0);
\draw[over](2,-1.4) arc (270:450:1.6) arc(270:90:0.1) arc(450:270:1.8) arc(270:90:0.1);

\draw[over] (1.5,1.6)arc (360:270:0.2) --+(-0.4,0);
\draw[over] (1.7,1.6)arc (360:270:0.4) --+(-0.4,0);
\draw[over] (1.5,-1.2)arc (0:90:0.2) --+(-0.4,0);
\draw[over] (1.7,-1.2)arc (0:90:0.4) --+(-0.4,0);
\draw[over] (2,1.9)arc (0:90:0.1) --+(-0.2,0);
\draw[over] (0.9,1.9)arc (0:90:0.1) arc(90:100:1.8);
\draw[over] (0.9,-1.5)arc (360:270:0.1) arc(270:260:1.8);
\draw[over] (1.9,1.9)arc (180:270:0.1) arc(450:360:1.6);
\draw[over] (2,-1.5)arc (360:270:0.1) --+(-0.2,0);
\end{scope}

\begin{scope}[yshift=-10.5cm]

\draw[over] (-0.2,0.2)arc(180:270:1)--+(0.4,0);
\draw[over] (-0.4,0.2)arc(180:270:1.2);
\draw[over] (-0.8,0.2)arc(180:260:1.6);
\draw[over] (-1,0.2)arc(180:260:1.8);
\end{scope}
\end{tikzpicture}
\caption{An example of a GBR: the two-fold Bing double of the Hopf link.}
\label{figure:GBR}
\end{figure}

Write $X_L := S^3 \sm N(L)$ for the exterior of $L$.
Let $\mu_{i} \subset X_L$ be an oriented meridian of the $i$th component of $L$, and let $\lambda_i \subset X_L$ be a zero-framed oriented longitude.  Both are smoothly embedded curves.  Make $\mu_i$ small enough that $\lk(\mu_i,\lambda_i)=0$  (of course $\lk(\mu_i,L_i) =1$ and $\lk(\lambda_i,L_i)=0$).
Let $N(\mu_i),\, N(\lambda_i) \subset X_L$ be closed tubular neighbourhoods, each homeomorphic to $S^1 \times D^2$.

A \emph{Round handle} $H$ is a copy of $S^1 \times D^2 \times D^1$.  The \emph{attaching region} is $S^1 \times D^2 \times S^0 \subset \partial(S^1 \times D^2 \times D^1) \cong S^1 \times S^2$.  The notion of round handles is due to Asimov~\cite{Asimov}. 

\begin{definition}
  Given an $m$-component link $L$, construct a manifold $R(L)$ by attaching $m$ round handles $\{H_i\}_{i=1}^m$ to $D^4$ as follows.  For the $i$th round handle, glue $S^1 \times D^2 \times \{-1\}$ to $N(\mu_i) \subset X_L \subset S^3= \partial D^4$, and glue $S^1 \times D^2 \times \{1\}$ to $N(\lambda_i)$.  In both cases use the zero-framing for the identification of $N(\mu_i)$ and $N(\lambda_i)$ with $S^1 \times D^2$.
  Note that the link $L$ lies in $\partial R(L)$.
\end{definition}

\noindent The key question will be whether $L$ is slice in $R(L)$.

\begin{definition}[Round Handle Slice]
  A link $L$ is \emph{Round Handle Slice $($RHS$)$} if $L \subset \partial R(L)$ is slice in $R(L)$, that is if $L$ is the boundary of a disjoint union of locally flat embedded discs in $R(L)$.
\end{definition}

\begin{theorem}\label{theorem:RHP}
Suppose that the topological surgery and $s$-cobordism conjectures hold for free fundamental groups.  Then for any link $L$ with pairwise linking numbers all zero, $L$ is round handle slice.
\end{theorem}

\begin{problem}
The Round Handle Problem is to determine whether all pairwise linking number zero links are round handle slice.
\end{problem}

By Theorem~\ref{theorem:RHP}, a negative answer for one such link would contradict the logical union of the topological surgery conjecture and the $s$-cobordism conjecture for free fundamental groups.
It is suggested by Freedman and Krushkal, but by no means compulsory, to focus on the links arising as GBRs.
It is also suggested that one might try to adapt Milnor's invariants to provide obstructions.
The primary purpose of this problem, like the AB slice problem, is to provide a way to get obstructions to surgery and $s$-cobordism.  Key work on the AB slice problem includes~\cite{Freedman-AB,Freedman-Lin:1989-1,Krushkal-AB,Freedman-Krushkal-2016-A}.

We briefly recall the statements of these conjectures and their relation to the disc embedding problem.

\begin{conjecture}[Topological surgery conjecture]\label{conjecture:surgery}
  Every degree one normal map $(M,\partial M) \to (X,\partial X)$ from a compact $4$-manifold $M$ to a $4$-dimensional Poincar\'{e} pair $(X,\partial X)$, that is a $\Z[\pi_1(X)]$-homology equivalence on the boundary, is topologically normally bordant rel.\ boundary to a homotopy equivalence if and only if the surgery obstruction in $L_4(\Z[\pi_1(X)])$ vanishes.
\end{conjecture}

\begin{conjecture}[$s$-cobordism conjecture]\label{conjecture:s-cobordism}
   Every compact topological $5$-dimensional $s$-cobordism $(W;M_0,M_1)$, that is a product on the boundary, is homeomorphic to a product $W \cong M_0 \times I \cong M_1 \times I$, extending the given product structure on the boundary.
\end{conjecture}

In Section~\ref{section:disc-embedding} we will explain why the union of these two conjectures is equivalent to the disc embedding conjecture, stated below.  In the statement of this conjecture we use the equivariant intersection form
\[\lambda \colon H_2(M,\partial M;\Z[\pi_1(M)]) \times H_2(M;\Z[\pi_1(M)]) \to \Z[\pi_1(M)]\]
  and the group-valued self-intersection number
  \[\mu \colon H_2(M;\Z[\pi_1(M)]) \to \frac{\Z[\pi_1(M)]}{g \sim w(g)g^{-1},\, 1 \sim 0},\]
  where $w \colon \pi_1(M) \to C_2 = \{\pm 1\}$ is the orientation character.  Also note that the transverse spheres are required to be \emph{framed}, which means that they have trivialised normal bundles.

\begin{conjecture}[Disc embedding conjecture]\label{conjecture:disc-embedding}
    Let $f_i\colon (D^2,S^1) \looparrowright (M,\partial M)$ be a collection of generically immersed discs in a compact $4$-manifold $M$ with disjointly embedded boundaries. Suppose that there are framed generically immersed spheres $g_i \colon S^2 \looparrowright M$ such that for every $i,j$ we have   $\lambda(g_i,g_j)=0$, $\mu(g_i)=0$, and the $g_i$ are transverse spheres, so $\lambda(f_i,g_j) = \delta_{ij}$.  Then the circles $f_i(S^1)$ bound disjointly embedded, locally flat discs in $M$ with geometrically transverse spheres, inducing the same framing on $f_i(S^1)$ as the $f_i$.
\end{conjecture}

Conjectures~\ref{conjecture:surgery}, \ref{conjecture:s-cobordism}, and \ref{conjecture:disc-embedding} are already theorems for good groups, a class of groups containing groups of subexponential growth~\cite{Freedman-Teichner:1995-1,Krushkal-Quinn:2000-1}, and closed under taking subgroups, extensions, quotients, and direct limits.

\begin{remark}
	 The obstruction theory presented in the proof of \cite[Lemma~5.4]{Freedman-Krushkal-2016-A}, which forms part of the proof given there of Theorem~\ref{theorem:RHP}, is incomplete.  First, $H^3(R,\partial R;\pi_2(R')) \cong H_1(R;\pi_2(R')) = 0$ since $\pi_2(R')$ is a free $\Z[\pi_1(R')] \cong \Z[\pi_1(R)]$-module, so the obstruction here certainly vanishes, as asserted in \cite{Freedman-Krushkal-2016-A}. However a potentially non-trivial obstruction, not considered in \cite{Freedman-Krushkal-2016-A}, lies in $H^4(R,\partial R;\pi_3(R'))$.  Analysing this depends on the relationship between the intersection forms of $R$ and $R'$. Our proof of Theorem~\ref{theorem:RHP} avoids obstruction theory altogether.
\end{remark}

\begin{remark}
  Our proof implies that every \emph{knot} is round handle slice, since in that case the proof applies Conjectures~\ref{conjecture:surgery} and~\ref{conjecture:s-cobordism} with fundamental group $\Z$. But for fundamental group $\Z$ these conjectures are theorems, since they are both implied by the disc embedding theorem~\cite[Section~2.9,~Theorem~5.1A]{Freedman-Quinn:1990-1}.
\end{remark}

\subsection*{Acknowledgements}

We thank Allison N.\ Miller for delivering a lecture on the Round Handle Problem that motivated us to produce this alternative proof. We also thank Peter Feller and Arunima Ray for helpful comments on a previous version, prompting us to elucidate the relation between the disc embedding conjecture and the topological surgery and $s$-cobordism conjectures. We also thank an anonymous referee for several useful comments that helped improve the exposition.

We are grateful to the Hausdorff Institute for Mathematics in Bonn, in whose
fantastic research environment this paper was written. The first author was partly supported by NRF grant 2019R1A3B2067839.
The second author was supported by an NSERC Discovery Grant.

\section{Proof of Theorem~\ref{theorem:RHP}}

The proof of Theorem~\ref{theorem:RHP} involves the construction of an $s$-cobordism rel.\ boundary from the manifold $R(L)$, henceforth abbreviated to $R$, to another $4$-manifold $R'$, in which $L$ is slice.
We begin with a Kirby diagram for $R$, shown in Figure~\ref{figure:RHP2}.

\begin{figure}[htb]
\centering
\begin{tikzpicture}
\begin{knot}[
clip width=10, clip radius=4pt,ignore endpoint intersections=false,]

\strand[thick] (3.2,0.3)
to [out=up, in=right] (1,2.5)
to [out=left, in=up] (-1.2,0.3)
to [out=down, in=left] (0.5,-2)
to [out=right, in=left] (1,0.3)
to [out=right, in=left] (1.5,-2)
to [out=right, in=down] (3.2,0.3);

\strand[thick]
(0.5,-1.4)
to [out=up,in=left] (1,-1.1)
to [out=right, in=up] (1.5,-1.4)
to [out=down, in=right] (1,-1.7)
to [out=left, in=down] (0.5,-1.4);

\strand[thick,color=red]
(1,-0.2)
to [out=right, in=down] (1.8,0.6)
to [out=up, in=right] (1,1.4)
to [out=left, in=up] (0.2,0.6)
to [out=down,in=left] (1,-0.2);

\strand[thick]
(1,1.2)
to [out=right, in=down] (1.35,2)
to [out=up, in=right] (1,2.8)
to [out=left, in=up] (0.65,2)
to [out=down, in=left] (1,1.2);

\draw[color=black] (1.4,3.1) node {$d$};
\draw[color=black] (-1.45,0.2) node {$0$};
\fill[color=black] (0.5,-1.4)  circle[radius=1.7pt];
\draw[color=red] (2,0.8) node {$L_i$};
\flipcrossings{2,4,6,7,9}
\end{knot}
\end{tikzpicture}
\caption{A handle diagram for $R$ in $N(L_i)$. Replicate this for each $i=1,\dots,m$.}
\label{figure:RHP2}
\end{figure}

First we will explain the figure, then we will explain why this is a diagram for $R$.  The diagram does not show the literal Kirby diagram for $R$.  Rather, the curve labelled $d$ specifies a solid torus, as the complement of a regular neighbourhood of this curve.  Inside the solid torus a dotted circle, corresponding to a 1-handle, and a zero-framed circle, corresponding to a $2$-handle, can be seen.  Embed a copy of this solid torus into a closed tubular neighbourhood $N(L_i)$ for each $i=1,\dots,m$, using the zero framing. One therefore has $m$ 1-handles and $m$ 2-handles, one pair in each solid torus neighbourhood $N(L_i)$, arranged as shown in Figure~\ref{figure:RHP2}.   The diagram also shows the link component $L_i$ parallel to the core of the solid torus.

Now we explain why Figure~\ref{figure:RHP2} is a diagram for the $4$-manifold $R$.
A round handle can be constructed from a 1-handle and a 2-handle whose boundary goes around one attaching circle of the round handle (a meridian of $L$), traverses the 1-handle, goes around the other attaching circle (a zero-framed longitude of the same component of $L$), and then traverses the $1$-handle in the other direction.   Ignoring the link $L$, we see that $R$ is diffeomorphic to the zero-trace of $L$ with $m$ $1$-handles added.

\begin{figure}[htb]
\centering
\begin{tikzpicture}
\begin{knot}[
clip width=10, clip radius=4pt,ignore endpoint intersections=false,]

\strand[thick] (3.2,0.3)
to [out=up, in=right] (1,2.5)
to [out=left, in=up] (-1.2,0.3)
to [out=down, in=left] (0.5,-2)
to [out=right, in=left] (1,0.3)
to [out=right, in=left] (1.5,-2)
to [out=right, in=down] (3.2,0.3);

\strand[thick]
(0.5,-1.4)
to [out=up,in=left] (1,-1.1)
to [out=right, in=up] (1.5,-1.4)
to [out=down, in=right] (1,-1.7)
to [out=left, in=down] (0.5,-1.4);

\strand[thick,color=red]
(1,-0.2)
to [out=right, in=down] (1.8,0.6)
to [out=up, in=right] (1,1.4)
to [out=left, in=up] (0.2,0.6)
to [out=down,in=left] (1,-0.2);

\strand[thick]
(1,-0.8)
to [out=right, in=down] (2.4,0.6)
to [out=up, in=right] (1,2)
to [out=left, in=up] (-0.4,0.6)
to [out=down,in=left] (1,-0.8);

\strand[thick]
(1,1.2)
to [out=right, in=down] (1.35,2)
to [out=up, in=right] (1,2.8)
to [out=left, in=up] (0.65,2)
to [out=down, in=left] (1,1.2);

\strand[thick,color=blue]
(2.1,1.3)
to [out=right, in=down] (2.3,1.5)
to [out=up, in=right] (2.1,1.7)
to [out=left, in=up] (1.9, 1.5)
to [out=down, in=left] (2.1,1.3);

\draw[color=black] (1.4,3.1) node {$d$};
\draw[color=black] (-0.6,0.8) node {$0$};
\draw[color=black] (-1.45,0.2) node {$0$};
\fill[color=blue] (2.26,1.6)  circle[radius=1.7pt];
\fill[color=black] (0.5,-1.4)  circle[radius=1.7pt];
\draw[color=red] (2,0.8) node {$L_i$};
\flipcrossings{2,4,6,8,9,11,13,15}
\end{knot}
\end{tikzpicture}
\caption{The handle diagram from Figure~\ref{figure:RHP2} with a cancelling pair introduced.}
\label{figure:RHP1}
\end{figure}

Figure~\ref{figure:RHP1} shows another diagram for $R$ with a cancelling 1-handle and 2-handle pair introduced in each $N(L_i)$.

Next, Figure~\ref{figure:RHP3} shows a Kirby diagram, with the same convention as above, for a 4-manifold that we call $R_M$.  Here $M$ stands for ``middle,'' since this manifold will lie in the middle of the $s$-cobordism we are about to construct.

\begin{figure}[htb]
\centering
\begin{tikzpicture}
\begin{knot}[
clip width=10, clip radius=4pt,ignore endpoint intersections=false,]

\strand[thick] (3.2,0.3)
to [out=up, in=right] (1,2.5)
to [out=left, in=up] (-1.2,0.3)
to [out=down, in=left] (0.5,-2)
to [out=right, in=left] (1,0.3)
to [out=right, in=left] (1.5,-2)
to [out=right, in=down] (3.2,0.3);

\strand[thick]
(0.5,-1.4)
to [out=up,in=left] (1,-1.1)
to [out=right, in=up] (1.5,-1.4)
to [out=down, in=right] (1,-1.7)
to [out=left, in=down] (0.5,-1.4);

\strand[thick,color=red]
(1,-0.2)
to [out=right, in=down] (1.8,0.6)
to [out=up, in=right] (1,1.4)
to [out=left, in=up] (0.2,0.6)
to [out=down,in=left] (1,-0.2);

\strand[thick]
(1,-0.8)
to [out=right, in=down] (2.4,0.6)
to [out=up, in=right] (1,2)
to [out=left, in=up] (-0.4,0.6)
to [out=down,in=left] (1,-0.8);

\strand[thick]
(1,1.2)
to [out=right, in=down] (1.35,2)
to [out=up, in=right] (1,2.8)
to [out=left, in=up] (0.65,2)
to [out=down, in=left] (1,1.2);

\strand[thick,color=blue]
(2.1,1.3)
to [out=right, in=down] (2.3,1.5)
to [out=up, in=right] (2.1,1.7)
to [out=left, in=up] (1.9, 1.5)
to [out=down, in=left] (2.1,1.3);

\draw[color=black] (1.4,3.1) node {$d$};
\draw[color=black] (-0.6,0.8) node {$0$};
\draw[color=black] (-1.45,0.2) node {$0$};
\draw[color=blue] (1.75,1.3) node{$\alpha_i$};
\draw[color=blue] (2.5,1.4) node{$0$};
\draw[color=black] (2.25,-0.55) node{$\beta_i$};
\draw[color=black] (-1.4,-0.7) node{$\gamma_i$};
\fill[color=black] (0.5,-1.4)  circle[radius=1.7pt];
\draw[color=red] (2,0.8) node {$L_i$};
\flipcrossings{2,4,6,8,9,11,13,15}
\end{knot}
\end{tikzpicture}
\caption{A handle diagram for $R_M$ in $N(L_i)$. The 2-handles in the picture are labelled $\alpha_i$, $\beta_i$ and $\gamma_i$.}
\label{figure:RHP3}
\end{figure}

The diagram for $R_M$ is very similar to the diagram for $R$ from Figure~\ref{figure:RHP1}; in order to get from the diagram for $R_M$ to that for $R$, inside each solid torus neighbourhood $N(L_i)$, change the zero-framed 2-handle whose attaching curve is labelled $\alpha_i$ in Figure~\ref{figure:RHP3} to a 1-handle. That is, for each $i$, perform surgery on the 2-sphere obtained from the core of the 2-handle union a disc bounded by the attaching circle in $D^4$. The fundamental group of $R_M$ is $\pi_1(R_M) \cong F_m$, the free group on~$m$ letters, generated by meridians of the dotted circles.

Note that, by virtue of the cores of the $\beta_i$ $2$-handles, $L$ is slice in $R_M$.  To see this, observe that $L_i$ can be passed through the attaching region of the $\alpha_i$ $2$-handle in Figure~\ref{figure:RHP3}.

\begin{remark}
  This remark  explains why the remainder of the proof is necessary: it is far from obvious that $L$ is slice in $R$.
In Figure~\ref{figure:RHP1}, $L_i$ cannot be passed through a dotted circle corresponding to a $1$-handle, so the argument just given cannot be used to show that $L_i$ is slice in $R$ via Figure~\ref{figure:RHP1}.  On the other hand if one isotopes the link through the attaching region of a $2$-handle, one cannot later use the core of that $2$-handle to construct an embedded slice disc, so one cannot use Figure~\ref{figure:RHP2} to see that~$L$ is slice in~$R$.
\end{remark}

Next, there are also generically immersed 2-spheres in $R_M$ obtained from the union of the cores of the $\beta_i$ 2-handles with immersed discs~$D_i$ in~$D^4$ bounded by the $\beta_i$ attaching curves.   By choosing the immersed discs $D_i$ so that their normal bundles induce the 0-framing on the curves $\beta_i\subset S^3$, we have framed immersed spheres. We call these the $\beta_i$-spheres.  The linking number zero hypothesis implies that the algebraic intersection numbers in $\Z[\pi_1(R_M)] \cong \Z[F_m]$ between these 2-spheres vanish.

Consider similar framed spheres arising from the round handle $2$-handles, namely the $2$-handles whose attaching curves are labelled $\gamma_i$ in Figure~\ref{figure:RHP3}.  For each $i$, isotope the curve $\gamma_i$ through $\beta_i$ and out from the 1-handle; that is, pull the oxbow part straight until $\gamma_i$ is a round circle, parallel to $\beta_i$.  Embed the isotopies in a collar $S^3 \times [1-\varepsilon,1]$.   Use parallel push offs of the discs $D_i$, minus their intersection with $S^3 \times [1-\varepsilon,1]$, to cap the resulting curves. We have just constructed discs $E_i$ with boundary $\gamma_i$, that intersect the discs $D_j$ algebraically in $\delta_{ij}$.  Cap off the discs $E_i$ with the cores of the $\gamma_i$ 2-handles to obtain framed immersed 2-spheres in $R_M$, that we call the $\gamma_i$-spheres.
The $\beta_i$- and $\gamma_j$-spheres are algebraically dual over $\Z[F_m]$.

\begin{lemma}\label{lemma:emb-spheres}
There exist framed, locally flat, embedded spheres $B_i \subset R_M$ in the complement of the slice discs for $L$, with $B_i$ regularly homotopic to the $\beta_i$-sphere for $i=1,\dots,m$.
\end{lemma}
\begin{proof}
To prove Lemma~\ref{lemma:emb-spheres}, we will apply the disc embedding conjecture to immersed Whitney discs~$f_k$ pairing up double points of the $\beta_i$-spheres, in the complement of the slice discs for $L$ in $R_M$, and in the complement of the $\beta_i$-spheres themselves. We will then perform the Whitney move using the resulting embedded Whitney discs to obtain the spheres $B_i$.

We argue that the immersed Whitney discs $f_k$ can be found.
First, apply the geometric Casson lemma~\cite[Lemma~3.1]{Freedman:1982-1},~\cite[Section~1.5]{Freedman-Quinn:1990-1} to convert the $\beta_i$-spheres and the $\gamma_j$-spheres from algebraic duals into geometric duals, intersecting in precisely one point if $i=j$ and with empty intersection otherwise.

Preliminary immersed Whitney discs $f_k'$ can be found in the complement of slice discs for $L$ because the slice discs for $L$ in $R_M$ use push offs of the core of the $\beta_i$ 2-handles, whereas the double points of the $\beta_i$-spheres lie in the interior of $D^4$. So one can find immersed Whitney discs in $D^4$ pairing up all double points among the $\beta_i$-spheres.  However, these initial Whitney discs $f_k'$, which we can assume to be framed Whitney discs by boundary twisting, might intersect the $\beta_i$-spheres.  Tube each intersection of a Whitney disc with a $\beta_i$-sphere into a parallel copy of the dual sphere $\gamma_i$.  This produces Whitney discs $f_k$ in $R_M$ that are framed and disjoint from both the slice discs for $L$ and the $\beta_i$-spheres.

Construct framed transverse spheres for the $f_k$ from Clifford tori for the double points, with caps given by normal discs to the $\beta_i$-spheres tubed into the dual $\gamma_i$-spheres. Use the caps to symmetrically contract~\cite[Section~2.3]{Freedman-Quinn:1990-1} the tori to immersed spheres. See \cite[Corollary~5.2B]{Freedman-Quinn:1990-1} for more details. Call the resulting spheres $g_k$.
All intersections among the transverse spheres $g_k$ arose from contraction, so they cancel algebraically over $\Z[F_m]$, and we therefore have $\lambda(g_k,g_{\ell})=0 = \mu(g_k)$ for every $k,\ell$.  Similarly, all of the intersection points between the $f_k$ and the $g_{\ell}$ cancel, except those arising from the original intersection points between Clifford tori and the Whitney discs $f_k$. It follows that the $f_k$ and the $g_{\ell}$ are algebraically dual over $\Z[F_m]$. We may therefore apply the disc embedding Conjecture~\ref{conjecture:disc-embedding} to find embedded Whitney discs, in the complement of the slice discs for $L$ and in the complement of the~$\beta_i$-spheres. The disc embedding conjecture has no hypothesis on the fundamental group, so we do not need to control the fundamental group here.  Whitney moves across the embedded discs resulting from Conjecture~\ref{conjecture:disc-embedding} give a regular homotopy to the desired framed embedded spheres~$B_i$.
This completes the proof of Lemma~\ref{lemma:emb-spheres}.
\end{proof}

Perform surgery on $R_M$ using these framed embedded spheres $B_i$, and define $R'$ to be the $4$-manifold obtained as result of these surgeries.
Note that $L$ is still slice in $R'$, since the spheres $B_i$ lie in the complement of the slice discs.

\begin{lemma}\label{lemma:s-cob}
The $4$-manifolds $R$ and $R'$ are $s$-cobordant rel.\ boundary.
\end{lemma}

\begin{proof}
To prove Lemma~\ref{lemma:s-cob}, start with $R_M$.  The trace of surgeries on the $\alpha_i$-spheres gives a cobordism to $R$.  The trace of surgeries on the $\beta_i$-spheres gives a cobordism to $R'$.  The union of the two cobordisms along $R_M$ is an $s$-cobordism from $R$ to $R'$, since algebraically the intersection numbers $\alpha_i \cdot \beta_j = \delta_{ij}$.  This completes the proof of Lemma~\ref{lemma:s-cob}.
\end{proof}

Note that we used duals to the $\beta_i$-spheres twice, once to apply surgery and once to prove that we have an $s$-cobordism. However we use \emph{different} duals.  For the surgery we used the $\gamma_i$-spheres arising from the round handle $2$-handles.  For the $s$-cobordism, we used the $\alpha_i$-spheres.

Then since $R$ and $R'$ are $s$-cobordant, the $s$-cobordism Conjecture~\ref{conjecture:s-cobordism} implies that they are homeomorphic rel.\ boundary.  Since the homeomorphism is an identity on the boundary, the link~$L$ is preserved.  Thus the image of the slice discs for $L$ in $R'$ under the homeomorphism $f \colon R' \to R$ are slice discs for $L$ in~$R$.  It follows that $L$ is Round Handle Slice as desired.
This completes the proof of Theorem~\ref{theorem:RHP}.

\section{Disc embedding is equivalent to surgery and $s$-cobordism}\label{section:disc-embedding}

In this section we briefly argue that the disc embedding Conjecture~\ref{conjecture:disc-embedding} is equivalent to the combination of the surgery and $s$-cobordism Conjectures, numbered \ref{conjecture:surgery} and \ref{conjecture:s-cobordism} respectively. There are no new equivalences described in this section. Indeed, references are given throughout, mostly to the relevant subsections of \cite{Freedman-Quinn:1990-1}.  We include this section for readers wanting a succinct guide to establishing these equivalences.

 We will argue that the following are equivalent: (i) surgery and $s$ cobordism; (ii) disc embedding; (iii) height 1.5 capped gropes contain embedded discs with the same boundary; (iv) certain links $L \cup m$, to be described below, are slice with standard slice discs for $L$.
 We will show:
 \[\text{(i)} \underset{(4)}{\implies} \text{(iv)} \underset{(3)}{\iff} \text{(iii)} \underset{(2)}{\iff} \text{(ii)} \underset{(1)}{\implies} \text{(i)}.\]


\begin{enumerate}[(1)]
\item \emph{The disc embedding conjecture \textup{(ii)} implies \textup{(i)} surgery and $s$-cobordism.} This follows from inspection of the high dimensional proof: the proof of topological surgery in dimension four and the five dimensional topological $s$-cobordism theorem can be reduced to precisely the need to find embedded discs with geometrically transverse spheres in the presence of algebraically transverse spheres. See for example~\cite{Luck-surgery-book} for an exposition of the high dimensional theory.  The $s$-cobordism theorem requires an extra argument to find the transverse spheres, which can be found in~\cite[Chapter~7]{Freedman-Quinn:1990-1}.
\item \emph{The disc embedding conjecture \textup{(ii)} is equivalent to the statement \textup{(iii)} that every height $1.5$ capped grope contains an embedded disc with the same framed boundary.}  For one direction, if disc embedding holds, then we can use it to find a disc in a height $1.5$ capped grope, as follows.  The caps on the height $1$ side are immersed discs, and parallel copies of the symmetric contraction of the height $1.5$ side, together with annuli in neighbourhoods of the boundary circles, give transverse spheres that have the right algebraic intersection data. See \cite[Section~2.6]{Freedman-Quinn:1990-1} for the construction of transverse gropes within a grope neighbourhood, which are then symmetrically contracted~\cite[Section~2.3]{Freedman-Quinn:1990-1} to yield transverse spheres.  Apply disc embedding to find embedded discs with framed boundary the same as the height 1 caps' framed boundary. These correctly framed embedded discs can be used to asymmetrically contract the first stage of the height $1.5$ grope to an embedded disc.  On the other hand, a collection of discs with transverse spheres as in Conjecture~\ref{conjecture:disc-embedding} gives rise to a height $1.5$ capped grope with the same boundary and with geometrically transverse spheres for the bottom stage, as shown in \cite[Section~5.1]{Freedman-Quinn:1990-1}.  Thus if every height $1.5$ capped grope contains an embedded disc, then disc embedding holds.
\item \emph{Height 1.5 capped gropes contain embedded discs with the same boundary \textup{(iii)} if and only if \textup{(iv)} certain links $L \cup m$ are slice with standard slice discs for $L$.}
A Kirby diagram for a capped grope consists of an unlink $L$, in the form of a link obtained from the unknot by iterated ramified Bing doubling, followed by a single operation of ramified Whitehead doubling. Place a dot on every component to denote that they correspond to $1$-handles; a neighbourhood of a capped grope is diffeomorphic to a boundary connected sum of copies of $S^1 \times D^3$.  The boundary circle of the grope is represented by a meridian $m$ to the original unknot.  One can think of performing the ramified Bing and Whitehead doubling on one component of the Hopf link.  A grope contains an embedded disc with the same framed boundary if and only if this link $L \cup m$ is slice with standard smooth slice discs for all the dotted components. The desired embedded disc is the slice disc for~$m$. See \cite[Proposition~12.3A]{Freedman-Quinn:1990-1} for further details.
\item \emph{Surgery and $s$-cobordism \textup{(i)} together imply \textup{(iv)} that the links $L \cup m$ are slice with standard slice discs for $L$.}  Let $L \cup m$ be any link from the family constructed in the previous item, using iterated ramified Bing and Whitehead doubling on one component of the Hopf link.  The zero surgery on $L \cup m$ bounds a spin $4$-manifold over a wedge of circles since the Arf invariants of the components vanish.  By the topological surgery conjecture, this can be improved, via a normal bordism rel.\ boundary, to be homotopy equivalent to the  wedge of circles.   Attach a $2$-handle to fill in the surgery torus $D^2 \times S^1$ of~$m$.  The remaining 4-manifold is homeomorphic to a boundary connected sum of copies of $S^1 \times D^3$, by the $s$-cobordism conjecture. Therefore it is homeomorphic to the exterior of standard smooth slice discs for $L$ in $D^4$. (We have no control over the remaining slice disc, whose boundary is the link component~$m$.)  Thus surgery and $s$-cobordism imply that the link $L \cup m$ is slice with standard slice discs for $L$.  More details are given in \cite[Section~11.7C]{Freedman-Quinn:1990-1} and the preceding sections of Chapter $11$.
\end{enumerate}

\bibliographystyle{alpha}
\def\MR#1{}
\bibliography{research}
\end{document}